\documentclass[smallextended]{svjour3}

\usepackage[utf8x]{inputenc}

\usepackage{dsfont}

\usepackage{amsmath,amssymb}
\usepackage{pst-all,pstricks-add}
\usepackage{color}
\usepackage{graphicx}

\usepackage{placeins}
\usepackage{listings}
\usepackage{subfigure}

\title{Mutant number distribution in an exponentially growing population}
\author{Peter Keller \and Tibor Antal}

\institute{
School of Mathematics, University of Edinburgh, Edinburgh, EH9 3FD, UK
}

\newcommand{\IP}{\mathds{P}}
\newcommand{\IE}{\mathds{E}}
\newcommand{\Fze}[3]{\mathop{F\/}\nolimits\!\left(\genfrac{}{}{0pt}{0}{#1}{#2};#3\right)}

\newcommand{\LPSM}{\lim_{\genfrac{}{}{0pt}{}{{N\rightarrow\infty,\mu\rightarrow0}}{\theta\text{ const.}}}}

\newcounter{romcount}

\newcommand{\genfunc}[2]{g_{#1}(#2)}

\newcommand{\loggen}[2]{\Lambda_{#1}(#2)}
\newcommand{\charexp}[2]{\mathcal E_{#1}(#2)}
\DeclareMathOperator{\sgn}{sgn}
\DeclareMathOperator{\EW}{\mathds E}
\DeclareMathOperator{\var}{Var}
\date{\today}

\begin{document}

\maketitle

\begin{abstract}
We present an explicit solution to a classic model of cell-population growth introduced by Luria and Delbr\"uck
\cite{LD1943} 70 years ago to study the emergence of mutations in bacterial populations. In this model a wild-type
population is assumed to grow exponentially in a deterministic fashion. Proportional to the wild-type population size,
mutants arrive randomly and initiate new sub-populations of mutants that grow stochastically according to a
supercritical birth and death process. We give an exact expression for the generating function of the total number of
mutants at a given wild-type population size. We present a simple expression for the probability of finding no mutants,
and a recursion formula for the probability of finding a given number of mutants.
In the ``large population-small mutation''-limit we recover recent results of Kessler and Levin \cite{kessler2013} for a
fully stochastic version of the process. 
\end{abstract}

\section{Introduction}
When a population of bacteria is attacked by a lethal virus, often a sub-population survives. At the beginning of the
1940's an important question was whether this resistance is due to adaptation which is induced
under the stress of attack, or is simply due to mutations that occurred beforehand during the expansion of the
population. To clarify this question, Luria and Delbrück conducted their now famous experiments in 1943 \cite{LD1943},
and showed that indeed the natural variability of cells can withhold a sub-population from extinction. They
formulated a simple mathematical model in which both wild-type and the mutant cells grow deterministically, but the
mutants appear randomly, proportional to the wild-type population size. They derived many properties of
the model, in particular for the mutant size distribution and proposed a method to estimate the mutation rate from data.

In the seminal paper of Lea and Coulson \cite{LC1949}, the original model was extended to allow stochastic growth
of the mutant population as a pure birth process. They derived the distribution of the number of mutants for the first
time for neutral mutation. In 1955 Bailey published elegant computations and some results on his own
modifications of the process, see \cite{Bailey1955}. Since then many efforts have been undertaken to understand the
process better. The review paper of Zheng \cite{zheng1999}, gives a formidable overview of the history of
the process and clarifies most concerns related to the infinite moments of the proposed distributions. 

New interest has kindled recently in the mutant distribution of the fully stochastic version, where wild-type cells
grows according to a birth and death process. Including cell death into the model extended the range of its possible
applications. This model was formulated by Kendall \cite{kendall60}, and a full solution was provided in
\cite{AntalKrap2011}, where the Kolmogorov equations for the generating function of both cell types were solved
explicitly. From the generating function the joint probability of a given number of wild type and mutant cells can be
obtained for any finite times. Expressions for finite times become important for experimental studies
\cite{Clayton:2007aa,antal10,bozic2013}, where the asymptotic limit might be out of reach.

In many situations, most notably in the study of mutations in tumor growth \cite{vogelstein13,nowak06}, the age of the 
wild-type population is rarely known. At tumor detection we have a fairly good idea about the size of the tumor but the
time of its initiation is unknown. This led to studies of the mutant distribution at a fixed size of wild-type
population.  For neutral mutataion and pure birth processes this problem was solved by Angerer \cite{angerer2001}. 
Iwasa, Nowak and Michor \cite{iwasa2006} extended this model to non-neutral mutants and to birth-and-death
processes. They derive mutant distributions and resistance probability assuming the product of population size and
mutation rate to be small. Komarova suggested a very elegant method to obtain an approximate mutant distribution
\cite{komarova07}. 
More recently, in two remarkable papers \cite{kessler2013,KessleronArXiV} Kessler and Levin obtained the full mutant
distribution for a large but fixed size wild-type population. They used approximate methods to simplify the Kolmogorov
equations, and in an independent derivation they also used the exact solution of the fully stochastic case given in
\cite{AntalKrap2011}.
By letting the previously constant product of mutation rate and population size go to infinity, they derive
\(\alpha\)-stable distributions. For the same limit, similar results were derived with other methods
by Durrett and Moseley \cite{durrettmoseley} for beneficial mutations. This result was already given, but not proven, by
Mandelbrot in
\cite{mandelbrot1974}. Moehle treats the classic
case of neutral mutation utilizing Compound-Poisson-Processes in \cite{moehle2005}. In \cite{janson} Janson treats a
similar model with fixed, non-random number of offspring, by mapping a reducible multi-type branching process to
P\'olya-Urns and investigates several limits. 

In this paper we make the assumption that the wild-type population grows according to a deterministic
exponential function. Leaning on the formalism for arbitrary growth functions introduced in 
\cite{dewanji2005} and reviewed in Section 2. In Section 3 we rewrite the general integral representation of the
generating function of the number of mutants explicitly in terms of hypergeometric functions. We consider the special
case of neutral mutations separately in Section 4. After getting rid of the integral representation, we investigate
limits of the mutant distribution in section 5. Indeed, we recover all corresponding results from the above mentioned
papers for general parameters and extend them to finite wild-type populations, to mutants with explicit death, and to
deleterious (disadvantageous compared to wild-type) mutations. We give a recursion to calculate the probability
distribution of the mutants efficiently and analyze the distribution's tail behavior in-depth in the final section,
thereby extending the results
of \cite{pakes1993,prodinger1996}.

\section{General population size functions}

We consider a cell population that consists of two types of cells, a wild type (type A) and a mutant (type B).
Each $A$-cell independently of all other cells produces a mutant $B$-cell at rate $\nu$. If we approximate the size of
the $A$-cell population via the deterministic function \(f(t)\), so that mutants are produced at rate $\nu f(t)$, the
arrival times of new mutants follow a non-homogeneous Poisson process. Each \(B\)-cell descended from an \(A\)-cell
at time \(s<t\) is the initiator of a new sub-population of mutants (a clone), whose size we denote by \(Y_k\). 
At time \(t\) the total number \(K\) of clones is a Poisson random variable with mean
\begin{equation}
     m 
      = \EW(K)
      = \int_0^t \nu f(s) ds.
\label{eq:general mean}
\end{equation}
We assume that clones develop indepently as some stochastic process with generating function \(\genfunc
tz=E\left(z^{Y}\right)\).
Since each clone \(Y_k\) is generated according to a Poisson-Process, the family \((Y_i)_{i\in\{1,\ldots,K\}}\) is
independent, identically distributed (iid) and the
generating function of each clone is 
\begin{equation}
 \psi (z) 
   = \EW\left(z^Y\right)
   = \frac{\nu}{m} \int_0^t f(s) \genfunc{t-s}z ds.
\label{eq:psi}
\end{equation}
The total number $B_t$ of mutants at time $t$ is a Compound Poisson random variable
\begin{equation*}
 B_t 
   = \sum_{i=1}^K Y_i.
\end{equation*}
Using conditional expectation, the generating function of $B_t$ can be written as
\begin{equation*}
 G(z) 
   = \EW\left(z^{B_t}\right)
   = \EW\left(\EW(z^{B_t}|K)\right).
\end{equation*}
Now
\begin{equation*}
   \EW(z^{B_t}|K=k)
    = \psi^k(z)
\end{equation*}
since the clones are independent and thus 
\begin{equation}
 G(z) 
  = \EW(\psi^K)=\sum_{k\geq 0}\frac{(\psi(z)m)^k}{k!}e^{-m}=e^{m(\psi(z)-1)},
\end{equation}
which appears in \cite{moehle2005} and is characteristic for Compound Poisson variables. Using \eqref{eq:general mean}
and \eqref{eq:psi}, we can also write
\begin{equation}
G(z)=\exp\left(\nu\int_0^tf(s)\left[\genfunc{t-s}z-1\right]ds\right)
\label{Ggen}
\end{equation}
which appears in \cite{dewanji2005} in a more general setting.

Since the generating function of \(B_t\) is of exponential form, we introduce the following notation
for arbitrary random variable \(X\)
\begin{equation*}
 \loggen Xz=\log \EW(z^X)=\log G_X(z)
 \label{eq:definition log-gen function}
\end{equation*}
and refer to \(\loggen Xz\) as the \emph{log-generating function} of \(X\).

\section{Generating function for exponential growth}\label{sec:exponential growth}

Let us consider the special case of an exponentially growing wild-type population, such that $f(t)=e^{\delta  t}$, for
some \(\delta>0\).  
Hence mutants are produced at rate $\nu e^{\delta  t}$. Moreover, let us assume that each clone
behaves like a linear birth-death process with birth rate $\alpha $ and death rate $\beta $ with positive relative
fitness $\lambda =\alpha -\beta >0$, i.e.~the process is supercritical. The extinction
probability of a mutant clone is $q=\beta /\alpha = 1-\lambda /\alpha $, and its generating function is also
well known \cite{arthreyaney}
\begin{equation*}
 \genfunc sz = 1 - \frac{1-q}{1-\xi e^{-\lambda  s}}, \quad \xi=\frac{q-z}{1-z} = 1-\frac{1-q}{1-z}.
\end{equation*}
The pure birth case of \(\beta=0\) and thus \(q=0\) is well studied and corresponds to the assumption that cells only
divide, but never die.

We are interested in the distribution of the number of mutants at the time when the number of \(A\)-cells
reaches exactly \(N\). Since the \(A\)-cells grow deterministically, this happens at time 
$\tau=\log(N)/\delta$. We use the shorthand notation $B\equiv B_\tau$ for the number of mutants at time $\tau$.
Therefore the mutant log-generating function \eqref{Ggen} becomes
\begin{equation}
\begin{split}
\loggen{B}z 
   = \nu\int_0^\tau e^{\delta s} (\genfunc{\tau-s}z-1) ds
   = \frac{\mu }{\gamma }\int_{0}^\tau \frac1{e^{\lambda s}N^{-1/\gamma}\xi-1}\delta e^{\delta s}ds
\label{eq:dewanji finite}
\end{split}
\end{equation}
where 
\begin{equation*}
 \gamma=\delta/\lambda\quad \text{ and } \quad\mu=\nu/\alpha.
\end{equation*}
After a change of variable $u=e^{\delta s}/N$ we can compute the integral
\begin{equation}
\begin{split}
 \loggen{B}z 
& = \frac{-N\mu}{\gamma}\int_{1/N}^1\frac{1}{1-u^{1/\gamma}\xi}du\\
& = \frac{-N\mu}{\gamma}\int_{1/N}^1\sum_{k\geq0}(u^{1/\gamma}\xi)^k du
  = -N\mu\sum_{k\geq 0}\xi^k\left.\frac{u^{(k+1)/\gamma}}{\gamma+ k}\right|_{1/N}^1\\
& = \mu\sum_{k\geq0}\xi^k\left[\frac{N^{-k/\gamma}}{\gamma+k}-\frac{N}{\gamma+k}\right].
\end{split}
\label{eq:sumform}
\end{equation}
We can rewrite \eqref{eq:sumform} in terms of the hypergeometric function \eqref{eq:definition gauss
hypergeometric function}
\begin{equation}
\loggen{B}z 
 = \log G_{B}(z)
 = \frac{N\mu}{\gamma}
    \left[
      \frac1N\Fze{1,\gamma}{1+\gamma}{\xi N^{-1/\gamma}}-\Fze{1,\gamma}{1+\gamma}{\xi}
    \right],
\label{eq:generating function}
\end{equation}
where we utilized the Pochhammer symbol  $(\gamma)_k$ \eqref{eq:definition pochhammer symbol} to verify that
\begin{equation*}
 \frac{\gamma}{\gamma+k}=\frac{(\gamma)_k}{(1+\gamma)_k}.
\end{equation*}
The above equation \eqref{eq:generating function} is the final, exact, closed-form solution for the mutant distribution
for an exponentially growing wild type population. In the rest of the paper we shall analyze its properties.

The mean and variance of the number of mutants can be calculated by taking the usual approach of differentiating the
generating function \eqref{eq:generating function} or by using Dewanji's general expressions for mean and variance for
the case of arbitrary
growth function \(f(t)\), see \cite{dewanji2005}; this results in
\begin{equation}
 \EW (B) 
   = \frac{N\mu}{1-q}\cdot
   \begin{cases}
    \log N                         & \gamma=1\\
     \frac{1}{1-\gamma}(N^{1/\gamma-1}-1).  & \gamma\neq 1
   \end{cases}
 \label{eq:mean B_tau}
\end{equation}
and
\begin{equation}
 \var(B)
   = 
    \frac{N\mu}{(1-q)^2}\cdot
    \begin{cases}
       2(N-1)-(1+q)\log N                               & \gamma=1\\     
       (1+q)(N^{-1/2} -1)+\log N                            & \gamma=2\\
       \frac{2}{2-\gamma}N^{2/\gamma-1} +
       \frac{1+q}{\gamma-1}N^{1/\gamma-1}+
       \frac{q(2-\gamma)+\gamma}{(2-\gamma)(1-\gamma)}      & \gamma\not\in\{1,2\}.
    \end{cases}
\label{eq:variance B_tau}
\end{equation}
These expressions generalize those given in Zheng \cite{zheng1999}
(replace $N=\exp(\delta  t)$ and $\delta \equiv\beta_1 $, $\lambda \equiv\beta_2$ in \cite[(52),(53)]{zheng1999}).
 
\begin{figure}[ht]
 \begin{center}
  \includegraphics[width=0.5\textwidth]{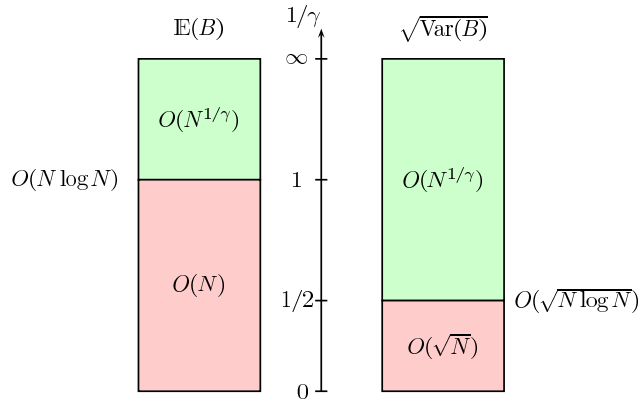}
 \end{center}
\caption{Orders of the mean number of mutants \(\EW(B)\) and its standard deviation \(\sqrt{\var(B)}\) for large 
$N$. Mutants have a fitness advantage for \(1/\gamma>1\) and a disadvantage for \(1/\gamma<1\) with respect to the
wild-type cells.}
\label{fig:overview orders}
\end{figure}

We give an overview of the large $N$ behavior of the expectation and the variance of the number of mutants $B$ in
figure
\ref{fig:overview orders}. The mean number of deleterious mutants  ($1/\gamma<1$)  is of the same order as the
wild type cells. However, the number of advantageous mutants ($1/\gamma>1$) is growing faster than the wild type
population. Note also that for advantageous mutants ($1/\gamma>1$) the mean and the standard deviation have the same
order, which implies that the fluctuations are important and a stochastic description of the process is essential. Only
for very deleterious mutants ($1/\gamma<1/2$) the process becomes self-averaging, and the fluctuations become as
predicted by the central limit theorem. For intermediate deleterious mutants ($1/2<1/\gamma<1$) the relative
standard deviation varies continuously with $\gamma$.

We can obtain the probabilities  $p_k=P(B=k)$ by Taylor expanding $G(z)$ in $z$, or by using the Gauss inversion
formula. Since this is computationally intense, we give a recursive formula for the probabilities instead
\begin{equation}
p_n
 = \begin{cases}
   e^{q_0}, &n=0\\
   \frac{1}{n} \sum_{k=0}^{n-1} (n-k) q_{n-k} p_k, & n\ge 1,
\end{cases}
\label{eq:recursion_general}
\end{equation}
where
\begin{equation}
 q_0 = \frac{N\mu}{\gamma} 
      \left[
          \frac1N\Fze{1,\gamma}{1+\gamma}{N^{-1/\gamma} q}
        - \Fze{1,\gamma}{1+\gamma}{q}
      \right]
\end{equation}
and for $k\ge 1$
\begin{equation}
\begin{split}
q_k 
  & =  \mu \sum_{j=1}^k \binom{k-1}{j-1} \frac{1}{j+\gamma}
       \left(
         \frac{1-q}{q-N^{1/\gamma}} 
       \right)^j
     \Fze{1,\gamma}{1+\gamma+j}{N^{-1/\gamma}q}\\
     &\quad + N\mu \frac{(k-1)!}{(\gamma+1)_{k}} \Fze{k,\gamma}{1+\gamma+k}{q}.\\
\end{split}
\label{general_rec_coeffs}
\end{equation}
We give a proof of this recursion in Appendix \ref{sec:rec}.

\section{Special case of neutral mutations, $\gamma=1$}

Often there is interest in mutations which do not change the behavior of the cell, so called \emph{neutral mutations}.
In this special case when $\gamma=1$, that is $\delta =\lambda$, we can further simplify the log-generating function
\(\loggen{B}z\) given in \eqref{eq:generating function}, by using \eqref{eq:log and hypergeo}. Alternatively, by using
the series expansion
\begin{equation*}
 \log(1-z)=-\sum_{k\geq 1}\frac{z^k}{k},
\end{equation*}
we can rewrite \eqref{eq:sumform} for $\gamma=1$ as 
\begin{equation}
\begin{split}
 \loggen{B}z 
   & = \mu\sum_{k\geq0}\xi^k\left[\frac{N^{-k}}{k+1}-\frac{N}{k+1}\right]
     =\frac{N\mu}{\xi}\log\left(\frac{1-\xi}{1-\xi/N}\right).
\label{eq:cgf for b=1}
 \end{split}
\end{equation}

Hence the generating function becomes
\begin{equation}
 G(z) 
   = \left(
        \frac{1-\xi}{1-\xi/N}
      \right)^{\frac{N\mu}{\xi}}.
\label{eq:generalized exact Lea-Coulson pgf}
\end{equation}
By introducing the variables
\begin{equation}
y=\frac{z-q}{1-q} = \frac{\xi}{\xi-1}, \quad \phi=1-\frac{1}{N}, \quad \theta= N\mu,
\label{eq:def y phi theta}
\end{equation}
we obtain the form
\begin{equation}
 G(z) 
   = (1-\phi y)^{\theta \frac{1-y}{y}}.
 \label{b_one_gener}
\end{equation}
If we further specialize to mutant cells that cannot die, that is $\beta=0$ (which implies $q=0$ and $y=z$), we find
\begin{equation}
 G(z) 
   = (1-\phi z)^{\theta \frac{1-z}{z}}.
 \label{b_one_nodeath}
\end{equation}
This formula was first derived in [Lea-Coulson], and also given in \cite{zheng1999} with some historical
perspective.

The coefficients in the recursion formula \eqref{eq:recursion_general} become simpler for $\gamma=1$ 
\begin{equation}
\label{eq:gamma1_rec}
 q_k = 
 \begin{cases}
 -\frac{\theta}{q} \log \left( 1+\phi \frac{q}{1-q} \right)& k=0\\
 \frac{\theta \phi^k}{[1-q(1-\phi)]^k} \left[ \frac{1}{k} - \frac{\phi}{k+1} \Fze{1,1}{2+k}{-\frac{\phi q}{1-q}} \right]
  & k\ge 1
 \end{cases}
\end{equation}
which is also derived in Appendix \ref{sec:rec}.
Note that 
\begin{equation}
 p_0 = P(B=0) = e^{q_0} =  \left( 1+\frac{q\phi}{1-q} \right)^{-\theta/q}
\end{equation}
and all other $p_n$ probabilities can be obtained from recursion \eqref{eq:recursion_general}.
When mutants do not die, i.e. $\beta=0$ (which implies $y=z$), the coefficients further simplify to
\begin{equation}
 q_k = 
 \begin{cases}
 -\theta\phi& k=0\\
 \theta \left( \frac{\phi^k}{k} - \frac{\phi^{k+1}}{k+1} \right)   & k\ge 0,
 \end{cases}
\end{equation}
as given in \cite[p.18]{zheng1999}.
Note that the coefficients \(q_k\) appear in \cite[Appendix D]{iwasa2006} as the approximated probabilities for large
\(N\) and \(N\mu\lll1\), in which case \(\phi\sim1\).

%

\section{Limit behavior}\label{sec:limit behaviour}

In Section \ref{sec:exponential growth} we derived a recursion for the exact probability distribution of the number of
mutants present at time \(\tau=\frac{\log N}{\delta}\), under the assumption of exponential wild-type-growth. The
coefficients of this recursion, given in \eqref{general_rec_coeffs}, are quite complicated, so we seek for an easier
limit
case. In \eqref{eq:generating function}, the first addend stems from the lower integral boundary \(0\) and should
vanish for large \(\tau\) resp.~large \(N\) and small mutation rate. This turns out to be true only if the
\(\theta=N\mu\) is held
constant. The  resulting distribution is a Compound Poisson random variable, enabling us to directly determine the
distribution of the size of a clone, i.e.~the sub-population size of mutants descended from one original mutant. This
distribution shows already a power-law tail,which will be discussed in Section \ref{sec:tail behaviour}. 

In applications, \(\theta\) can be large, hence we discuss a consecutive limit in which \(\theta\) goes to
infinity. Since a Compound Poisson random variable is the sum of a random number of i.i.d.\  variables, it is not
surprising that the distribution of the limit variable, which we call \(Z\), is \(\alpha\)-stable. Although the theory
of generalized limit theorems is rich, see for example \cite{durrettbook}, we use our explicit results for the
generating function to perform the limit directly. Since we consider only the convergence of
generating functions or Laplace-transforms, all limits in this paper are meant as convergence in distribution.

\begin{figure}[ht]
 \begin{center}
 \includegraphics[width=0.6\textwidth]{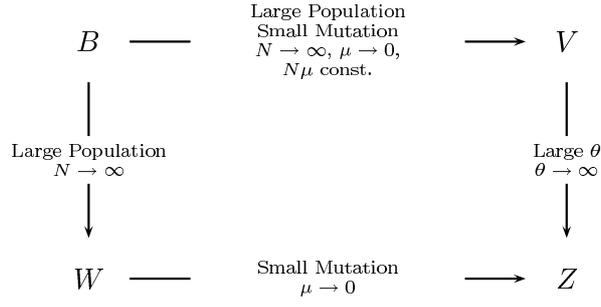}
 \end{center}
\caption{A schematic overview of the possible limits. All limits are convergence-in-distribution results.}
\label{fig:limits overview}
\end{figure}

Note that in \cite{durrettmoseley}, a non-rigorous proof was given for the direct limit from \(B\) to \(Y\) (see figure
\ref{fig:limits overview}) for the fully stochastic case in a slightly different setting for \(0<\gamma<1\).
It utilizes an approximation of the wild-type population size by a deterministic exponential growth function, that stems
from the fact that in the two-type, fully stochastic model the wild-type population behaves like a one-type birth and
death process if the mutation rate \(\nu\) is very small. Then the number of \(A\)-cells can be approximated
by \(\exp(\delta t)X\), \(X\sim Exp(1)\).  A similar limit in \cite{moehle2005} covers the \(\gamma=1\) case of the
classic Luria Delbr\"uck distribution. 

Another limit approach is to fix the mutation rate and let \(N\rightarrow\infty\) under 
a proper rescaling. Results for this limit were presented in \cite{mandelbrot1974} and their derivation given in a
privately distributed second part of the paper, which is no longer available.
We reproduce these results in terms of hypergeometric functions and then take the limit \(\mu\rightarrow 0\) under a
similar scaling as in the \(N\mu\)-constant case, to recover once again the limit variable \(Z\).

%

\subsection{Large Population-Small Mutation Limit}

In applications, the mutation rate \(\nu\) is usually small, while the population size \(N\) is very large. We
therefore investigate the simultaneous limit \(N\rightarrow \infty\) and \(\mu=\nu/\alpha\rightarrow0\), such that
\(\theta=N\mu\) is held constant. We call this limit \emph{Large Population-Small Mutation Limit} (LPSM). Note that
we introduced \(\theta\) already in \eqref{eq:def y phi theta}, in analogy to the notation introduced in
\cite{zheng1999}.

Formally we can express the LPSM-limit as 
\begin{equation*}
 \LPSM B=V,
\end{equation*}
where \(V\) is the limiting random variable, which we characterize via its log-generating function.
When \(\theta\) is held constant, the log-generating function of \(B\), given in \eqref{eq:generating function}, 
depends only in the first addend on \(N\). We therefore expand the first addend into a power series in
\(\xi\). For arbitrary \(\gamma>0\)
\begin{equation*}
 \begin{split}
  \frac 1N\Fze{1,\gamma}{1+\gamma}{N^{-1/\gamma}\xi}
    = \sum_{k\geq 0}\frac{\gamma}{\gamma+k}\xi^k\frac1{N^{k/\gamma+1}}
      \rightarrow 0,\text{ as \(N\rightarrow\infty\)}.
 \end{split}
\label{eq:log gen of X_theta}
\end{equation*}
Thus the immediate result for the LPSM-limit is 
\begin{equation}
 \loggen{V}z=\LPSM\loggen{B}z=-\frac\theta\gamma\Fze{1,\gamma}{1+\gamma}{\xi}.
 \label{eq:mv gen func}
\end{equation}
which we can rewrite in terms of
\(y=\frac{z-q}{1-q}\)
by using \eqref{eq:inversion hyperg 02}
\begin{equation}
 \loggen Vz=\frac{\theta}{\gamma}(1-y)\Fze{1,1}{1+\gamma}{y}.
 \label{eq:loggen v in y}
\end{equation}
This expression is the second addend in the general formula \eqref{eq:generating function}, so
we can adapt the recursion \eqref{eq:recursion_general} for the probability distribution of \(V\) easily. This yields
the recursion coefficients
\begin{equation}
 q_k=\begin{cases}
     -\frac\theta\gamma\Fze{1,\gamma}{1+\gamma}{q}&, k=0\\
      \theta \frac{(k-1)!}{(\gamma+1)_{k}}
\Fze{k,\gamma}{1+\gamma+k}{q}&,k\geq1.
     \end{cases}
\label{eq:coefficients of Lambda Xc}
\end{equation}
Note that the coefficients \(q_k\) appear in \cite[(16)]{iwasa2006} as the approximated probabilities for large \(N\)
and \(N\mu\lll1\).

The expectation and variance of \(V\) can be derived as usual via derivatives of the generating
function. The computations are tedious, in particular because the convergence behavior
of the hypergeometric function depends on \(\gamma\), but not very interesting. We give the results in Table
\ref{tab:exp and var of V}. Note that they are consistent with  the application of the \(N\rightarrow\infty\),
\(\mu\rightarrow0\) limit directly to the mean and variance of \(B\) given in \eqref{eq:mean B_tau} and
\eqref{eq:variance B_tau}. Interestingly, the mean is finite only for \(\gamma>1\) and the variance only for
\(\gamma>2\).

\begin{table}
\begin{center} 
\begin{tabular}{c|c|c}
&\(\EW(V)\)
&\(\var(V)\)\\
\hline
\(\gamma>2\)
&\(\frac{\theta}{(1-q)(\gamma-1)}\)
&\(\frac{\theta}{(1-q)^2}
      \left(\frac{q(2-\gamma)+\gamma}{(\gamma-2)(\gamma-1)}\right)\)\\
\(1<\gamma\leq2\)
&\(\frac{\theta}{(1-q)(\gamma-1)}\)
& \(\infty\)\\
\(0<\gamma\leq1\)
&\(\infty\)
&\(\infty\)    
 \end{tabular}
\end{center}
\caption{Overview of mean and variance of \(V\), depending on \(\gamma\).}
\label{tab:exp and var of V}
\end{table}

We note that the LPSM-limit is independent of the initial number \(N_0\) of wild-type cells. This can be easily
seen, when we choose \(f(t)=N_0 e^{\delta t}\) in \eqref{Ggen}. Then the integral representation of the
log-generating function \eqref{eq:dewanji finite}  is just \(G_{B}^{N_0}(z)\).  This property directly mimics the
branching property of a fully stochastic two-type branching process. Adapting the calculations from \eqref{eq:sumform},
the log-generating function
\begin{equation}
\loggen{B}z 
 = \log G^{N_0}_{B}(z)
 = \frac{N\mu}{\gamma}
    \left[
      \frac{N_0}N\Fze{1,\gamma}{1+\gamma}{\xi \frac{N}{N_0}^{-1/\gamma}}-\Fze{1,\gamma}{1+\gamma}{\xi}
    \right].
\label{eq:generating function with N0 wild types initially}
\end{equation}
Again, in the LPSM-limit only the first addend depends on \(N\) resp.~\(N_0\) and vanishes with \(N\rightarrow0\).

Note further, that the limit log-generating function  \(\loggen{V}{z}\) has a direct interpretation with respect to our
initial
model. We write
\begin{equation*}
  \loggen{V}{z}
    = -\frac\theta\gamma\Fze{1,\gamma}{1+\gamma}\xi
    = \frac\theta\gamma\int_0^1\frac{1}{u^{1/\gamma}\xi-1}du
\end{equation*}
by \eqref{eq:hypergeo as integral}. A change of variables with \(s=\frac{\log (uN)}{\delta}\) gives 
\begin{equation*}
 \frac\theta\gamma\int_0^1\frac{1}{u^{1/\gamma}\xi-1}du
    = \frac\theta{\gamma N}\int_{-\infty}^\tau\frac{1}{e^{\lambda s}N^{-1/\gamma}\xi-1}\delta e^{\delta s}ds.
\end{equation*}
Indeed, \(\loggen Vz\) is the generating function of the model started at \(-\infty\) instead of zero.

By another change of variables \(t=\tau-s\)
\begin{equation*}
 \loggen Vz 
   = \frac\theta\gamma\int_0^\infty\frac{1}{e^{-\lambda t}\xi-1}\delta e^{-\delta t}dt
   = \frac{\theta}{(1-q)\gamma}\int_0^\infty (\genfunc tz-1)\delta e^{-\delta t}dt.
\end{equation*}
Let \(X\sim Exp(\delta)\) be an exponential random variable with mean \(1/\delta\), then 
\begin{equation*}
 \loggen Vz=\frac{\theta}{(1-q)\gamma}\IE\left[\genfunc Xz-1\right].
\end{equation*}
This is the log-generating function of a Compound-Poisson random variable, where
\begin{equation*}
 \IE\left[\genfunc Xz\right]=1-(1-q)\Fze{1,\gamma}{1+\gamma}{\xi}.
\end{equation*}
Via the usual interpretation of a Compound Poisson, the limiting number of clones (i.e. the number of actual
mutation events) has a Poisson distribution with intensity \(\frac{\theta}{(1-q)\gamma}\) and \(\IE[\phi_X(z)]\)
describes the clone size distribution, that is the size of the population founded by a single original mutant.

\begin{figure}
\subfigure{\includegraphics[width=0.5\textwidth]{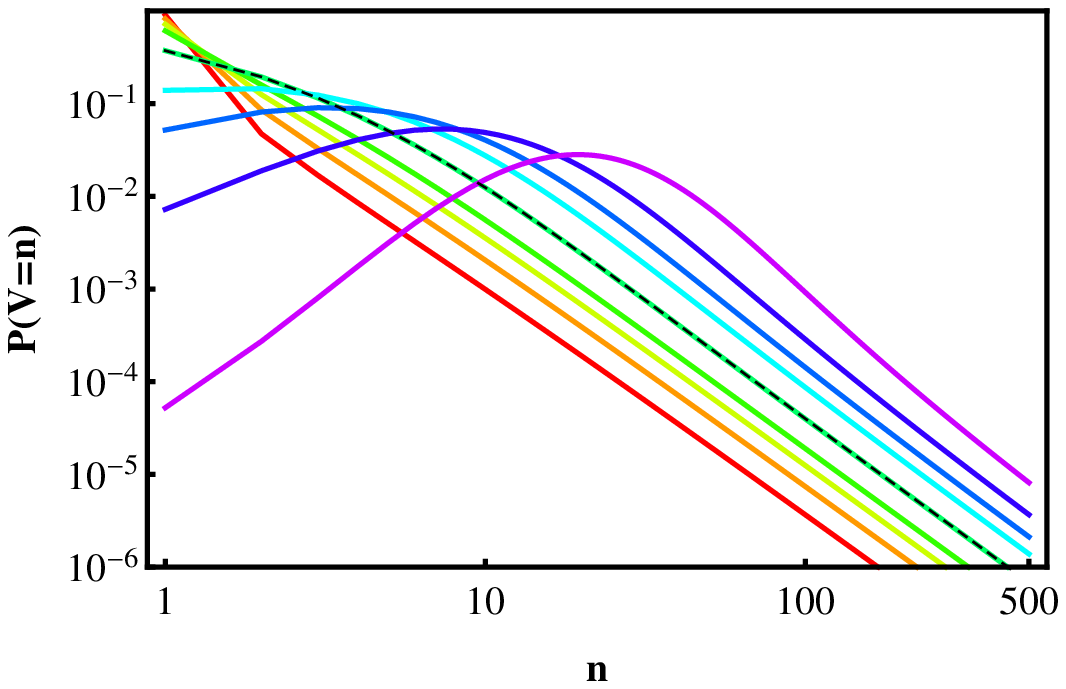}}
\subfigure{\includegraphics[width=0.5\textwidth]{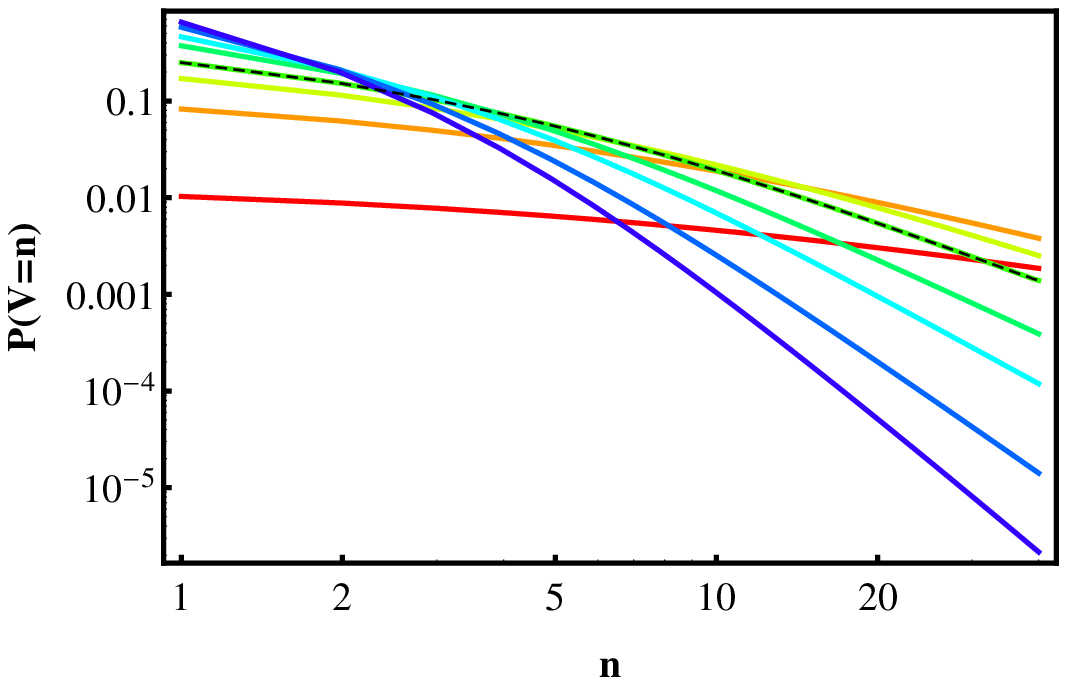}}

\caption{On the left, we plot the Probability distribution of the number of mutants in the large population-small
mutation limit with $\theta=N\mu$ constant. While keeping  \(\gamma=1.5\) and \(q=1/2\) constant, we vary
$\theta=$\(1/10\), \(1/5\), \(1/3\), \(1/2\), \(1\), \(2\), \(3\), \(5\), \(10\), beginning with the lowest (red) line.
On the right, we vary $\gamma$ choosing $0.25$, $0.5$, $0.75$, $1$, $1.5$, $2$, $3$ and $4$, beginning with the
horizontal (red) line, for \(\theta=1\) and \(q=1/2\). 
The case \(\gamma=1\) is indicated with a dashed line in both figures. We used recursion \eqref{eq:coefficients of
Lambda Xc} to calculate the probabilities.}
\label{fig:vary gamma theta}
\end{figure}

\begin{figure}
\begin{center}\includegraphics[width=0.35\textwidth]{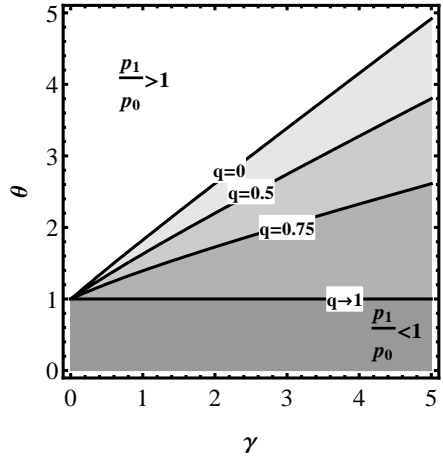}\end{center}
\caption{Comparison of the probability of no mutants \(p_0\) and a single mutant \(p_1\). Solid lines indicate
\(p_0=p_1\) depending on \(q\). We included the extreme case \(q=1\) as well.}
\label{fig:phase diagram}
\end{figure}

\begin{figure}
\begin{center}\includegraphics[width=0.4\textwidth]{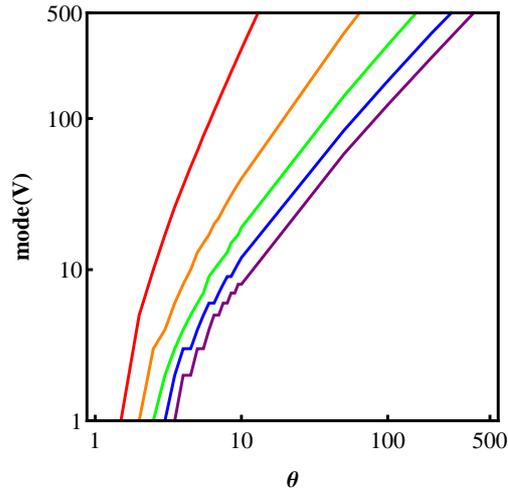}\end{center}
\caption{The mode of \(V\) plotted as function of \(\theta\) and \(q=0.5\). We chose (beginning from the left)
\(\gamma\) equal to \(0.5\), \(1\), \(1.5\), \(2\) and \(2.5\). The jagged appearance of some of the lines is due to
the mode taking only integer values.}
\label{fig:most probable mutant number constant gamma}
\end{figure}

\begin{figure}
\begin{center}\includegraphics[width=0.4\textwidth]{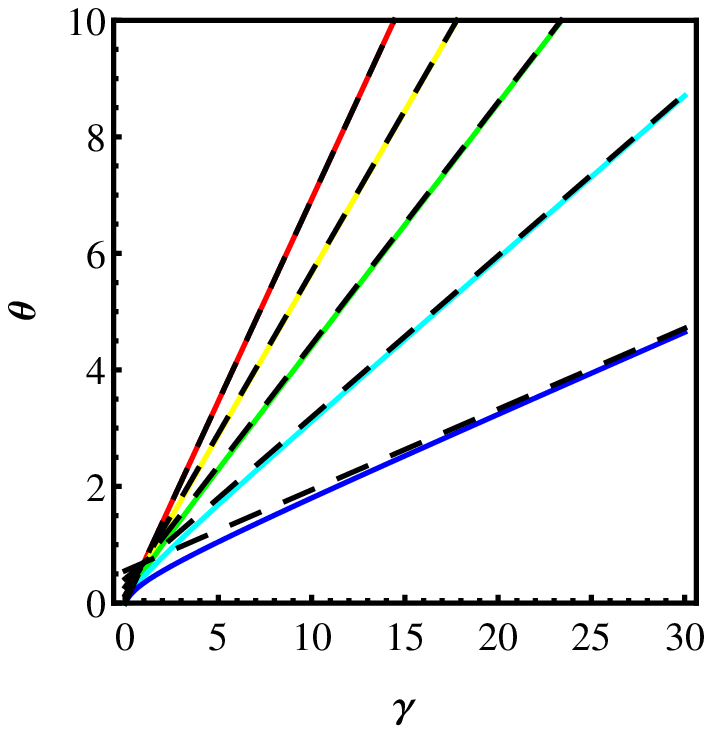} \end{center}
\caption{Plot of the contours dependend of \(\gamma\) and \(\theta\) where the probability \(p_0=1/2\) for
\(q=0,0.2,0.4,0.6,0.8\) (left to right) and the approximation (dashed line) given in \eqref{eq:niveau line p0}. Note
that the quality of the approximation depends on \(q\).}
\label{fig:resistance probability}
\end{figure}

We identified the distribution of \(V\) as (discrete) Compound Poisson random variable and argued that it
has infinite mean for \(0\leq\gamma\leq1\) and infinite variance for \(0\leq\gamma\leq2\), see Table \ref{tab:exp and
var of V}. Since the potentially non-finiteness of the moments makes the application of moment-methods very difficult
and unreliable, we give some information about the mode of \(V\). The mode is defined as the value k at which the
probability mass function of \(V\) takes its maximum. In uni-modal distributions the mode corresponds to the peak of
the
distribution.

From numerical analysis we are confident that the distribution of \(V\) is indeed uni-modal, as can be seen in
Figure \ref{fig:vary gamma theta} where we vary the parameters \(\gamma\) and \(\theta\). In general, however, it is
difficult to derive an explicit formula since we only have the recursion \eqref{eq:coefficients of Lambda Xc} available.
The recursion, however, makes the first few probabilities explicitly accessible. We therefore investigate the ratio
\begin{equation}
 \frac{p_1}{p_0}=\frac{\theta}{\gamma+1}\Fze{1,\gamma}{2+\gamma}{q}.
 \label{eq:ration p1/p0}
\end{equation}
This ratio is less than one, if the probability \(p_0\) is the maximum of the distribution and larger if the maximum is
bigger or equal to one. This ratio depends not only on \(\theta\) and \(\gamma\), but also on \(q\). In Figure
\ref{fig:phase diagram} we show a phase diagram on the \(\theta\)-\(\gamma\) plane displaying the boundary of regions
where no mutants are the most probable, for a few values of \(q\). The plot shows that the boundaries behave like a
linear function. This is can be confirmed, if we set \(p_1/p_0=1\) in \eqref{eq:ration p1/p0}, then \(\theta\) and
\(\gamma\) can be seperated and the boundary is described by
\begin{equation}
 \theta=(1+\gamma)\Fze{1,\gamma}{2+\gamma}{q}^{-1}\sim1+q+(1-q)\gamma+O(\gamma^2),\text{ for large \(\gamma\)}.
\label{eq:contourapproximation}
\end{equation}
The last approximation is due to an asymptotic result for the hypergeometric function, which we derive in
Appendix \ref{app:resistance}.
The ratio \(p_1/p_0\), however, gives no information about the actual value of
the mode(B), that is the most probable number of mutants. By numerically sampling the mode, using an implementation of
the recursion \eqref{eq:coefficients of Lambda
Xc} in Mathematica, we find that the mode increases rapidly for large \(\theta\) and small \(\gamma\), see also Figure
\ref{fig:most probable mutant number constant gamma}.

The probability \(p_0\) is indeed an important quantity in itself, since the ``resistance probability''
\(\IP(V>0)=1-p_0\) indicates, in the sense of the original Luria-Delbrueck formulation, the probability that a
population can escape extinction under the attack of a lethal virus due to the existence of resistent mutants. Thanks to
the recursion \eqref{eq:coefficients of Lambda Xc}, we
can give \(p_0\) explicitly 
\begin{equation}
 \IP(V=0)
     = \exp\left[-\frac{\theta}{\gamma}\Fze{1,\gamma}{1+\gamma}q\right].
\label{eq:resistance probability}
\end{equation}
The expression \(1-p_0\), derived differently, appears in \cite[(7)]{iwasa2006} (Note that they have another definition
of the mutation rate, which maps into ours with \(\delta=(1-u)r\) and \(\nu=ur\), where \(u\) is the mutation
probability). 
With increasing mode the resistance probability increases, in fact exponentially fast by \eqref{eq:resistance
probability}. For \(\gamma\rightarrow\infty\) and fixed \(\theta\) the
resistance probability goes to \(0\). For \(\gamma\) fixed and \(\theta\rightarrow\infty\), the probability goes to
one. Using a Taylor expansion of the hypergeometric function in \eqref{eq:resistance probability}, we can approximate
the resistance probability for large \(\gamma\) as
\begin{equation}
 -\log p_0\sim \frac{\theta}{\gamma(1-q)+q}.
\end{equation}
For a fixed \(p_0\), the two variables \(\theta\) and \(\gamma\) can be seperated and the contour associated to the
fixed value of \(p_0\) can be parameterized as 
\begin{equation*}
 \theta\sim-\left[\gamma(1-q)+q\right]\log p_0.
 \label{eq:niveau line p0}
\end{equation*}
A plot of the contours, where the resistance probability of \(V\) is equal to \(1/2\) for different \(q\) is given in
Figure \ref{fig:resistance probability} together with their approximations.

%

\subsection{Large values of $\mathbf\theta$}

We now let \(\theta\rightarrow\infty\), using the results 
of the previous section.
Obviously, for \(\theta\rightarrow\infty\) \eqref{eq:mv gen func} does not converge, thus we need to 
introduce a scaling. We formalize this convergence in distribution by
\begin{equation*}
 \lim_{\theta\rightarrow\infty}\frac{V}{a}-b=Z
\end{equation*}
where \(a,b\) depend only on \(\theta\). The scaling factor is
in fact already predetermined, see \cite[section 3.7, p.135]{durrettbook}, and is proportional to
\(\theta^{1/\min(\gamma,2)}\), but we will re-derive it for our case. We take the opportunity to scale out the survival
probability \(1-q\) wherever possible.
\newcommand{\expsc}{e^{-s/a}}

Let us abbreviate the logarithm of the Laplace transform of \(V\) by 
\begin{equation*}
\charexp{V}{s}=\log\EW(\exp(- V s))=\loggen V{e^{-s}}.
\label{eq:def of log laplace}
\end{equation*}
This definition directly implies
\begin{equation*}
\charexp{\frac{V}{a}-b}s=\charexp{V}{s/a}+bs.
\end{equation*}
We use this notation also for \(W\).

For \(y|_{z=\expsc}\), eq. \eqref{eq:loggen v in y} turns into
\begin{equation}
\charexp{V}{s/a}=\frac{\theta}{\gamma}\frac{\expsc-1}{1-q}\Fze{1,1}{1+\gamma}{\frac{\expsc-q}{1-q}}.
\label{eq:mv gen func simpl}
\end{equation}
We then apply equation \eqref{eq:inversion hyperg 01} to derive
\begin{equation*}
 \begin{split}
 &\charexp{V}{s/a}
     = \frac{\theta(\expsc -1)}{(1-q)(\gamma-1)}\Fze{1,1}{2-\gamma}{\frac{1-\expsc }{1-q}}
     +\frac{\pi}{\sin(\pi\gamma)}\left(-\frac{\theta^{1/\gamma}(\expsc-1)}{\expsc -q}\right)^\gamma
 \end{split}
\end{equation*}
For \(a\) to infinity \(\expsc-q\rightarrow 1-q\) and
\(\Fze{1,1}{2-\gamma}{\frac{1-\expsc}{1-q}}\rightarrow\Fze{1,1}{2-\gamma}{0}=1\) for all \(\gamma>0\).
With a Taylor expansion of \(\expsc\) around \(s=0\) we can write 

\begin{equation}
 \begin{split}
 \charexp{V}{s/a}
    & \sim \sum_{k\geq1}\frac{(-s)^k}{\theta^{-1}a^k(1-q)(\gamma-1) k!}
     +\frac{\pi}{\sin(\pi\gamma)}\left(\sum_{k\geq1}\frac{(-s)^k}{\theta^{-1/\gamma}a^k(1-q)
k!}\right)^\gamma\\
 \end{split}
\label{eq:asymptotic form}
\end{equation}

For the computations we consider four cases, which depend on the choice of \(\gamma\). The results are listed
in Table \ref{fig:c->oo limit results}. Note, that although the expressions for \(\gamma\in(0,1)\) and
\(\gamma\in(1,2)\) are equal, the means are not. Indeed this can be understood via a derivative with respect to \(s\)
and a limit \(s\rightarrow0\).


\begin{table}
 \begin{center}
  \begin{tabular}{c|c|c|c}
  \(\gamma\)
& \(\charexp{Z}{s}\)
& \(\EW(Z)\)
& \(\var(Z)\)\\
\hline 
  \(\gamma\in(0,1)\) 
& \(\frac{\pi}{\sin(\pi\gamma)}s^\gamma \)
& \(+\infty\) 
& \(+\infty\)\\
  \(\gamma=1\)       
& \(s\log s\)
& \(+\infty\) 
& \(+\infty\)\\
  \(\gamma\in(1,2)\)
& \(\frac{\pi}{\sin(\pi\gamma)}s^\gamma\)
& \(0\) 
& \(+\infty\)\\
  \(\gamma\geq 2\)
& \(\frac12s^2\)
& \(0\)
& \(1\)
\end{tabular}
 \end{center}
\caption{Overview or the results of the large \(\theta\) limit. Note that although the expressions in the first and
third row are identical, the means are not.}
\label{fig:c->oo limit results}
\end{table}

\colorbox{lightgray}{\(\mathbf{\gamma>2}\)}
For \(\gamma>2\) it is now intuitive to set \(a=\sqrt{\frac{\theta}{(\gamma-1)(1-q)}}\), since for this choice all terms
except the linear and quadratic term in the first addend of \eqref{eq:asymptotic form} vanish for
\(\theta\rightarrow\infty\). The linear term however diverges, so we compensate it by setting
\(b=\sqrt{\frac{\theta}{(1-q)(\gamma-1)}}\). Thus
\begin{equation*}
 \charexp{Z}{s}=\lim_{\theta\rightarrow \infty}\charexp{V}{s/a}+bs=\frac{s^2}{2},
\end{equation*}
which proofs that the limit random variable \(Z\) has a standard normal distribution.

\colorbox{lightgray}{\(\mathbf{0<\gamma<2, \gamma\neq1}\)}
In this case it is sufficient to set \(a=\theta^{1/\gamma}/(1-q)\), then \eqref{eq:asymptotic form} is
\begin{equation*}
 \charexp{V}{s/a}
    \sim -\frac{\theta^{1-1/\gamma}}{\gamma-1}s+O(\theta^{1-2/\gamma}) 
      +\frac{\pi}{\sin(\pi\gamma)}\left[s+O(\theta^{-1/\gamma})\right]^\gamma.
\end{equation*}
Therefore, if we set \(b=\frac{\theta^{1-1/\gamma}}{\gamma-1}\),
\begin{equation*}
 \charexp{Z}{s}=\lim_{\theta\rightarrow \infty}\charexp{V}{s/a}+bs=\frac{\pi}{\sin(\pi\gamma)}s^\gamma.
\end{equation*}
Note that for \(0<\gamma<1\), we could have chosen \(b=0\), since \(\frac{\theta^{1-1/\gamma}}{\gamma-1}\rightarrow0\).
The result resembles and in fact generalizes the result of \cite[Th. 6]{durrettmoseley} and \cite{mandelbrot1974}.

\colorbox{lightgray}{\(\mathbf{\gamma=1}\)}
Here we have to consider \eqref{eq:mv gen func simpl}, since the expansion \eqref{eq:asymptotic
form} has a singularity at \(\gamma=1\). We use  \eqref{eq:log and hypergeo} instead and gain
\begin{equation*}
\begin{split}
 \charexp{V}{s/a}
   & = \frac{\theta(\expsc-1)}{\expsc-q}\left[\log\left(\frac{\theta(\expsc-1)}{q-1}\right)-\log\theta \right].\\
\end{split}
\end{equation*}
With the same argumentation of Taylor expansion of \(\theta(\expsc-1)\) as in the proceeding cases we notice that
\(a=\theta/(1-q)\) is a sensible choice, since then
\begin{equation*}
 \charexp{V}{s/a}
   \sim s \log s-s\log\theta
\end{equation*}
and consequently for \(b=\log\theta\)
\begin{equation*}
 \charexp{Z}{s}=\lim_{\theta\rightarrow\infty}\charexp{V}{s/a}+bs=s\log s,
\end{equation*}
which also appears in \cite{moehle2005}.

\colorbox{lightgray}{\(\mathbf{\gamma=2}\)}
Here we face the same problem as in the previous case, since the expansion \eqref{eq:asymptotic form} is also singular
at
\(\gamma=2\). Via consecutive application of \eqref{eq:mv gen func simpl} and \eqref{eq:log and hypergeo 2} we gain
\begin{equation*}
 \charexp{V}{s/a}
   = \frac{\theta(\expsc-1)}{\expsc-a}-\frac{[\theta(\expsc-1)]^2}{\theta(\expsc-q)^2}
       \log\left(\frac{\theta\expsc-1}{\theta(q-1)}\right).
\end{equation*}
With the very same argumentation as in the \(\gamma=1\) case, we find for \(a=\frac{\sqrt{\theta\log\theta}}{1-q}\)
\begin{equation*}
 \charexp{V}{s/a}\sim -s\sqrt{\frac{\theta}{\log \theta}}+\frac{\log(\theta\log\theta)}{2\log\theta}s^2,
\end{equation*}
so that with \(b=\sqrt{\frac{\theta}{\log \theta}}\)
\begin{equation*}
 \charexp{Z}{s}=\lim_{\theta\rightarrow\infty}\charexp V{s/a}+bs=\frac12s^2.
\end{equation*}
Again, \(Z\) has a standard normal distribution.

\subsection{Large population limit at fixed mutation rate}

In this section we extend the results of Mandelbrot from \cite{mandelbrot1974} to allow the possibility of mutant
death. We consider a similar approach as in the Large \(\theta\)-limit but with mutation rate $\nu$ held constant, i.e.
\begin{equation*}
 \lim_{N\rightarrow\infty}\frac{B}{a}-b=W,
\end{equation*}
where \(a\) and \(b\) depend on \(N\).  
Like in \eqref{eq:def of log laplace}, we abbreviate with
\(\charexp Bs\) the log-Laplace transform of \(B\) and note that 
\begin{equation*}
 \charexp{\frac Ba-b}s=\charexp{B}{s/a}+ b s.
\end{equation*}
Applying first \eqref{eq:inversion hyperg 01} and then \eqref{eq:inversion hyperg 02} to \eqref{eq:generating function}
gives, noticing that \(\expsc-q\rightarrow 1-q\),
\begin{equation*}
 \begin{split}
  \charexp{B}{s/a}
     \sim \frac{\mu}{(\gamma-1)}
        \left[
           -\frac{N^{1/\gamma}(\expsc-1)}{1-q}
           \Fze{1,1-\gamma}{2-\gamma}{\frac{N^{1/\gamma}(\expsc-1)}{1-q}}
        \right.\\
        \left.
          + \frac{N(\expsc-1)}{1-q}\Fze{1,1-\gamma}{2-\gamma}{\frac{\expsc-1}{1-q}}
        \right].
 \end{split}
\end{equation*}
The parameter \(a\) controls the variance, therefore we know that we should choose
\(a\) of the order of \(\sqrt{\var(B)}\). The parameter \(b\) should be proportional to the order of the ratio
\(\EW(B)/\sqrt{\var(B)}\), see also Figure \ref{fig:overview orders}. The limits can now be computed like in the
LPSM-limit by expanding \(\expsc\) into a
power series. We give the chosen scaling factors together with the limit results and their mean and
variance in Table \ref{tab:limit results for the LPCM-limit}.
Note that we could have removed the mean for \(\gamma\in(0,1]\), but we need a non-zero expectation for the
\(\mu\rightarrow0\) limit, to recover \(Z\).

\begin{table}[ht]
\begin{center}
\begin{tabular}{c|c|c|c|c|c}
& \((1-q)a\)
& \(b\)
& \(\charexp Ws\)
& \(\EW(W)\)
& \(\var(W)\)\\
\hline
  $\gamma>2$
& $\sqrt N$
& $ \frac{\mu}{\gamma-1}\sqrt N$
& $\frac12\mu\frac{ \gamma-q(\gamma-2)}{(\gamma-2)(\gamma-1)}s^2$
& $0$
& $\frac{\mu(\gamma-q(\gamma-2))}{(\gamma-2)(\gamma-1)}$\\
  $\gamma=2$
& $\sqrt{N\log N}$
& $ \mu\sqrt{\frac{N}{\log N}}$
& $\mu\frac{s^2}2$
& $0$
& $\mu$\\
  $\gamma\in(1,2)$
& $N^{1/\gamma} $
& $\frac{\mu}{\gamma-1}(N^{1-\frac1\gamma}-1)$
& $\frac{\mu }{2-\gamma}s^2\Fze{1,2-\gamma}{3-\gamma}{-s}$
& $0$
& $\frac{2\mu}{2-\gamma}$\\
  $\gamma=1$
& $N$
& $\mu(1+\log N) $
& $\mu s(1+\log(1+s))$
& $\mu$
& $2\mu$\\
  $\gamma\in(0,1)$
& $N^{1/\gamma}$
& $0 $
& $\frac{\mu}{\gamma-1}s\Fze{1,1-\gamma}{2-\gamma}{-s} $
& $\frac{\mu}{\gamma-1} $
& $ \frac{2\mu}{2-\gamma}$\\
\end{tabular}

\end{center}
\caption{Overview of the limit results for \(N\rightarrow\infty\) and \(\mu=const\).}
\label{tab:limit results for the LPCM-limit}
\end{table}

In order to check that the two pathways for going from $B$ to $Z$ on Fig \ref{fig:limits overview} are equivalent, we
take the $\mu\to0$ limit of the above large $N$ limit. Again, we consider the limit
\begin{equation*}
 \lim_{\mu\rightarrow0}\frac{W}{a}=Z
\end{equation*}
where we assume for simplicity that \(a\) is now only dependent on \(\mu\). 

For \(\gamma\geq2\) there is nothing to show, since for any \(X\sim\mathcal N(\mu,\sigma^2)\) obviously
\(\frac{X-\mu}{\sigma}\sim \mathcal N(0,1)\). Therefore \(\frac{W-\EW(W)}{\var(W)}\) does not depend on \(\mu\) and the
limit for \(\mu\rightarrow0\) is trivially again a Gaussian.
We notice that by using \eqref{eq:simplification} in reverse
\begin{equation*}
 \frac{\mu}{2-\gamma}s^2\Fze{1,2-\gamma}{3-\gamma}{-s}
   =  \frac{\mu}{1-\gamma}s-\frac{\mu}{1-\gamma}s\Fze{1,1-\gamma}{2-\gamma}{-s}
\end{equation*}
and by eq. \eqref{eq:inversion hyperg 01} we get
\begin{equation*}
 \frac{\mu}{1-\gamma}s\Fze{1,2-\gamma}{3-\gamma}{-s}
  =\frac{\pi}{\sin(\pi\gamma)}\mu s^{\gamma} 
    - \frac{s\mu}{(s+1)\gamma}\Fze{1,1}{1+\gamma}{\frac{1}{1+s}}.
\end{equation*}
Finally, setting \(a=\mu^{1/\gamma}\), we get 
\begin{equation*}
 \charexp Zs=\lim_{\mu\rightarrow 0}\charexp{W}{s/a}=
  \begin{cases}
    s\log s&\gamma=1\\
    \frac{\pi}{\sin(\pi\gamma)}s^{\gamma} &\gamma\in(0,2)\setminus\{1\}
  \end{cases}
\end{equation*}
That shows that the scaling limit of \(W\) for \(\mu\rightarrow 0\) recovers indeed the distributions we got for the
Large \(\theta\) limit, which we give in Table \ref{fig:c->oo limit results}.

\subsection{On \(\alpha\)-stable distributions}

We point out that the limiting Laplace transforms given in Table \ref{fig:c->oo limit results} are well known
representatives of the class of \(\alpha\)-stable distributions, where \(\alpha\in(0,2]\). These distributions appear
as the limit distributions in the generalized version of the Central Limit Theorem, where iid. random variables are
added and rescaled, but the assumption of finite mean and variance is dropped. The Gaussian distribution is in this
context the extreme case for \(\alpha=2\) and the only limit distribution in this class with finite variance. As a
reference see \cite[sec. 2.7]{durrettbook}.

One of the many characterizations of \(\alpha\)-stable distributions is given via characteristic functions (the notation
varies strongly throughout the literature). A random variable \(X\) with characteristic function
\begin{equation*}
 \varphi_X(s)=\EW(\exp(iXs)).
\end{equation*}
is \(\alpha\)-stable resp.~\(X\sim S_\alpha(\sigma,\beta,\mu)\) if and only if
\begin{equation*}
 \log\varphi_X(s)
   = \begin{cases}
        -\sigma^\alpha|s|^\alpha(1-i\beta\sgn s\tan\frac{\alpha\pi}{2})+i\mu s   &\alpha\neq1\\
        -\sigma|s|(1+i\beta\frac{2}{\pi}\sgn s \log s)+i\mu s                   &\alpha=1.
     \end{cases}
\end{equation*}
By a formal substitution \(s\mapsto -is\), we can rewrite the distribution of \(Z\), 
as given in Table \ref{fig:c->oo limit results} using the above notation
\begin{equation}
\log \varphi_Z(s)
 =\begin{cases}
   -\sigma^\gamma
        {|s|}^\gamma\left(1-i \sgn(s)\tan\left(\frac{\pi\gamma}{2}\right)\right)& \gamma\in(0,2)\setminus\{1\}\\
   -\frac{\pi}{2}|s|\left(1+i\sgn(s)\frac{2}{\pi}\log s\right) &\gamma=1\\
   -\frac12s^2&\gamma\geq2,
  \end{cases}
\end{equation}
where 
\begin{equation*}
 \sigma=\left[\frac{\pi\gamma}{2}\csc\left(\frac{\pi\gamma}{2}\right)\right]^{1/\gamma}.
\end{equation*}
This means that \(Z\) is indeed an \(\alpha\)-stable random variable for all \(\gamma>0\), where the shape-parameter 
\(\alpha=\min(\gamma,2)\).

Unfortunately, the densities are unknown for the majority of \(\alpha\)-stable distributions. One exception is the
Gaussian
distribution (\(\alpha=2\)). The L\'evy-Distribution \(S_{1/2}(\sigma,1,0)\) with density
\begin{equation*}
 f(t)=\left(\frac{\sigma}{2\pi}\right)^{1/2}\frac{1}{t^{3/2}}\exp\left(-\frac{\sigma}{2t}\right)
\end{equation*}
and infinite moments and the Holtsmark distribution with finite mean are other cases (for the latter a lengthy
expression in terms of hypergeometric functions exists). 
In the special case for \(\gamma=1\), we find that \(Z\) is characterized by the
distribution \(S_1(\sigma,1,0)=S_1(\pi/2,1,0))\), which is the Landau-Distribution. This distribution was also
identified for the fully stochastic case by Kessler and Levin in \cite{kessler2013}.

\section{Tail behavior}\label{sec:tail behaviour}

The mutant distributions as we described them in the previous Section \ref{sec:limit behaviour} are quite complicated.
For this reason we study their tail behavior to gain some intuition. We use asymptotic
analysis methods as discussed in \cite{flajolet_book}.

We start with the most interesting case, the LPSM limit, where the number of mutants, denoted
by $V$ (see Figure \ref{fig:limits overview}), is characterized by generating function \eqref{eq:loggen v in y}.
The tail behavior
of $p_n=P(V=n)$ is encoded 
in $G(z)$ around its relevant (closest to origin) singularity.
This generating function is analytic at the origin, and its only singularity is at $z=1$, which implies
$y=\frac{z-q}{1-q}=1$, so we expand around that. We first use the transformation formula \eqref{eq:hyperg trafo
(1-z)} and then the
reflection formulas \eqref{eq:euler_reflection_B}, \eqref{eq:euler_reflection_C} to rewrite the log-generating function
as
$$
 \Lambda(z) = \theta \frac{1-y}{1-\gamma}  \Fze{1,1}{2-\gamma}{1-y}
 -\frac{\theta\pi}{\sin(\gamma\pi)} y^{-\gamma}(1-y)^\gamma.
$$
An expansion around $y=1$ gives
\begin{equation*}
\begin{split}
 \Lambda(z) &=  \frac{\theta}{1-\gamma} \left[ (1-y) + \frac{(1-y)^2}{2-\gamma} +O( (1-y)^3) \right]\\
 &- \frac{\theta\pi}{\sin(\gamma\pi)}
 \left[ (1-y)^\gamma + \gamma(1-y)^{1+\gamma}+O((1-y)^{2+\gamma})  \right],
\end{split}
\end{equation*}
and by using the series expansion of the exponential we get
\begin{equation*}
\begin{split}
 G(z) 
  & = e^{\Lambda(z)}
    = 1+\Lambda(z)+\frac{\Lambda^2(z)}{2} +O(\Lambda^3(z))\\
  & = K(1-y) - \kappa(1-y)^\gamma  - \kappa \left( \frac{\theta}{1-\gamma} +\gamma \right)
       (1-y)^{1+\gamma}  + \frac{1}{2}  \kappa^2 (1-y)^{2\gamma}\\
 &+ O((1-y)^{\{3\gamma,2+\gamma,1+2\gamma\})},
\end{split}
\end{equation*}
where $K$ is a polynomial with integer powers (hence do not contribute to the tail behavior), and for brevity we wrote
$$
 \kappa=\frac{\theta\pi}{\sin(\gamma\pi)}.
$$
Note that this is already an expansion in $1-z$ since
$$
 1-y=\frac{1-z}{1-q}.
$$
By using \eqref{eq:fla1} and \eqref{eq:euler_reflection_B} we obtain the large $n$ expansion
\begin{equation}
\label{eq:gengamma_q_series}
\begin{split}
 p_n &\sim  \frac{\theta\Gamma(1+\gamma)}{(1-q)^\gamma} n^{-1-\gamma}
 + \frac{\theta^2\kappa^2}{2\Gamma(-2\gamma)(1-q)^{2\gamma}}   n^{-1-2\gamma}\\
 &-\frac{\theta \Gamma(2+\gamma)}{(1-q)^{1+\gamma}} 
 \left( \frac{\theta}{1-\gamma} + \frac{\gamma(1+q)}{2} \right) n^{-2-\gamma}
 +O(n^{\{3\gamma,2+\gamma,1+2\gamma\}}).
\end{split}
\end{equation}
The leading order term has been derived in \cite{kessler2013}. The sub leading order is represented by the
$n^{-1-2\gamma}$ term for $\gamma<1$, and by the $n^{-2-\gamma}$ term for
$\gamma>1$.

Let us consider the above expansion at the special value \(\gamma=1\). 
As $\gamma\to 1$, both sub-leading terms in \eqref{eq:gengamma_q_series} diverge as $\propto1/|1-\gamma|$ but these two
singular terms cancel each other out, and one has to consider the next (constant) terms in the $(1-\gamma)$ series. A
logarithmic term appears through the expansion
$$
 n^{1-\gamma} = e^{(1-\gamma)\log n} = 1+(1-\gamma)\log n + O[(1-\gamma)\log n]^2
$$
for $\gamma\to1$.
We arrive at
\begin{equation}
\label{eq:gamma1death_final}
 p_n \sim \frac{\theta}{(1-q)n^2}
 + \frac{2\theta^2}{(1-q)^2} \frac{\log n}{n^3}
  + \frac{\theta^2[2C_E-3-2L(q)]-\theta(1+q)}{(1-q)^2 n^3}  + O\left( \frac{1}{n^4} \right)
\end{equation}
where $L(z)=\log\frac{1}{1-z}$.
For $q\to 0$, we have $L(q)\to 0$, and we obtain
\begin{equation}
\label{eq:asymp_simp}
\begin{split}
 p_n  &= \frac{\theta}{n^2} + 2\theta^2 \frac{\log n}{n^3} + \frac{\theta^2(2C_E-3)-\theta}{n^3}
+ 3\theta^3 \frac{\log^2 n}{n^4} + O\left(\frac{\log n}{n^4}\right)
\end{split}
\end{equation}
where $C_E=0.5772\dots$ is the Euler-Mascheroni constant. The first three terms have been calculated in
\cite{prodinger1996}.
Formally, all terms can be included, and here we just give the result
\begin{equation}
\begin{split}
 p_n &= \sum_{k=1}^\infty \frac{\theta^k}{k!(n+k)!} \lim_{r\to k} \partial_r^k \frac{\Gamma(n+k-r)}{\Gamma(-r)}\\ 
 &= \frac{\theta}{n(n+1)} + \frac{\theta^2[2C_E-3+2\Psi(n)]}{n(n+1)(n+2)} 
    + O\left( \frac{\theta^3 \log^2 n }{n^4} \right).\\
\end{split}
\end{equation}
This expression is exact for any specific order of $\theta$. By expanding it around $n=\infty$, using the expansion of
the Digamma function \(\psi(n)\) given in \eqref{eq:series_Psi}, we recover \eqref{eq:asymp_simp}.

The tail behavior of the full distribution of \(B\) with finite $N$ is similar, but additionally it
has an exponential cut-off. We demonstrate this behavior only for the $\gamma=1$ and $q=0$ case for simplicity.
In this case the generating function is given by \eqref{b_one_nodeath}, which is again
\begin{equation*}
 G(z) = (1-\phi z)^{\theta\frac{1-z}{z}}
\end{equation*}
with $\phi=1-1/N$. Now the singularity is at $z=1/\phi$, which leads to an exponential cut-off since
\begin{equation}
\label{eq:scaling_pole}
 [z^n] G(z) = \phi^n [z^n] G(z/\phi).
\end{equation}
where $[z^n]$ means the coefficient of $z^n$ in the Tailor expansion of the subsequent expression.
We now write
\begin{equation}
\label{eq:2factors}
 G(z/\phi) = (1-z)^{\theta\frac{\phi-z}{z}} = (1-z)^{-\mu} (1-z)^{(\theta-\mu)\frac{1-z}{z}}.
\end{equation}
The second factor is analytic at $z=0$, and its only singularity is at $z=1$ (and a branch cut from $z=1$ to infinity). 
Here the exponent approaches zero, so we expand first the exponential function
\begin{equation*}
\begin{split}
 G(z/\phi) &= \sum_{k\ge 0} \frac{(\mu-\theta)^k}{k!} (1-z)^{k-\mu} \left( \frac{L}{z} \right)^k\\
 &= (1-z)^{-\mu} +\left(\mu  - \theta \right) \frac{(1-z)^{1-\mu}}{z} L + \dots
\end{split} 
\end{equation*}
with $L\equiv L(z)=\log\frac{1}{1-z}$. Note that we didn't expand $1/z$ around $z=1$, and the reason for that becomes
clear in 
the next step. We obtain the coefficients by using \eqref{eq:fla1} and the Frobenius-Jungen method \eqref{eq:fla_bint}
from \cite{flajolet_book}
\begin{equation*}
\begin{split}
 [z^n] G(z/\phi) &= \sum_{k\ge 0} \frac{(\mu-\theta)^k}{k!} [z^{n+k}] (1-z)^{k-\mu} L^k\\
 &= \frac{n^{\mu-1}}{\Gamma(\mu)}\left( 1+\sum_{k\ge1} \frac{e_k}{n^k} \right)
 + \sum_{k\ge1} \frac{(\mu-\theta)^k}{k!} \frac{d^k}{da^k} \left. \binom{n+k+a-1}{n+k}\right|_{a=n-k}
 \end{split} 
\end{equation*}
The first few terms are
\begin{equation*}
\begin{split}
 [z^n] G(z/\phi) &= \frac{n^{\mu-1}}{\Gamma(\nu)} 
 \left[ 1+ \frac{\mu(\mu-1)}{2} n^{-1} + O(n^{-2}) \right]\\
 & +(\mu-\theta) \binom{n+\mu-1}{n+1} [\Psi(\mu-1+n)-\Psi(\mu-1)] + O(n^{\mu-3}\log n)
 \end{split} 
\end{equation*}
and the leading powers of $n$ are
\begin{equation*}
\begin{split}
 [z^n] G(z/\phi) &= \frac{n^{\mu-1}}{\Gamma(\mu)}
 + \frac{n^{\mu-2}}{\Gamma(\mu-1)} 
 \left\{  (\mu-\theta)[\log n -\Psi(\mu-1)]+\frac{\mu}{2} \right\}\\
 &+ O(n^{\mu-3}\log n).
 \end{split} 
\end{equation*}
Now we also include the exponential cut-off from \eqref{eq:scaling_pole} to obtain
\begin{equation}
\begin{split}
 p_n &\sim \frac{(1-1/N)^n}{\Gamma(\mu)}
  \left\{ \frac{1}{n^{1-\mu}} + 
  (1-\mu)(\theta-\mu) \frac{\log n }{n^{2-\mu}} \right. \\
   &- \left. (1-\mu)\left[(\theta-\mu)\Psi(\mu-1) +\frac{\mu}{2} \right] \frac{1}{n^{2-\mu}}
 + O\left( \frac{\log n}{n^{3-\mu}} \right) \right\}.
\end{split}
\end{equation}
The leading term of this expansion has been derived in \cite{pakes1993}.
Also, if we take the $\mu\to 0, N\to\infty$ limit of this asymptotic expansion, noting that $\Gamma(\mu)\sim1/\mu$ and
$\Psi(\mu-1)\sim-1/\mu$ as $\mu\to 0$, we recover the small mutation expression \eqref{eq:asymp_simp}.

\begin{figure}
\begin{center}
  \includegraphics[width=0.6\textwidth]{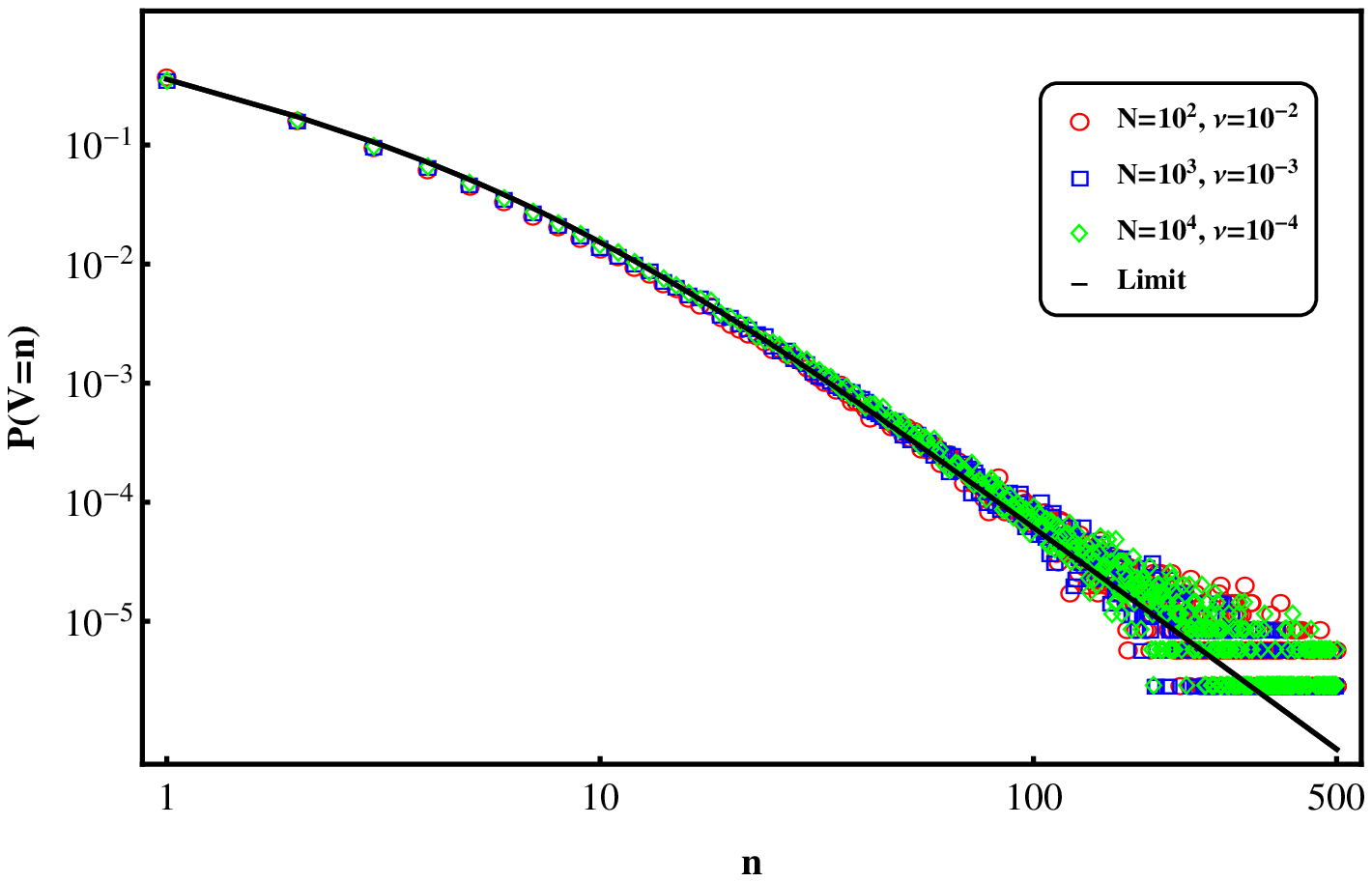}
\end{center}
\caption{Simulation of the fully stochastic model against the \(N\mu=1\)-limit for \(10^6\) trajectories each.}
\label{fig:Nmu1}
\end{figure}

\begin{figure}
\begin{center}\includegraphics[width=0.6\textwidth]{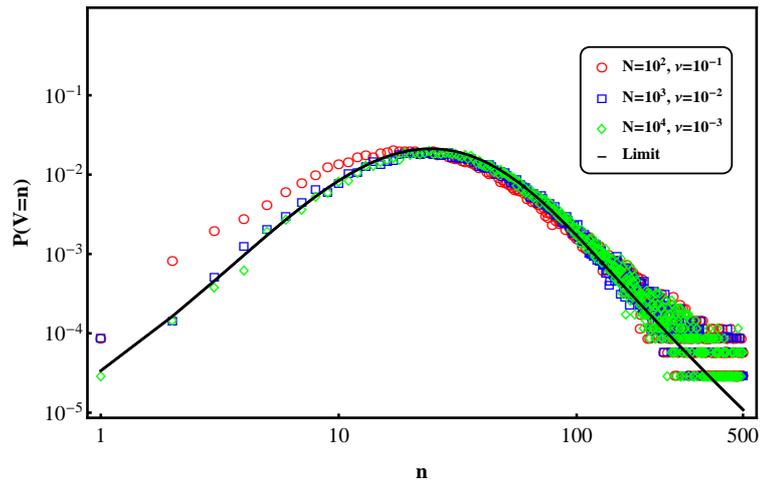}\end{center}
\caption{Simulation of the fully stochastic model against the \(N\mu=10\)-limit for \(10^6\) trajectories each.}
\label{fig:Nmu10}
\end{figure}

\section{Discussion and Summary}

We revisited a well-known semi-stochastic model with exponential wild-type growth and mutants that evolve according to a
supercritical birth and death process. We give an explicit expression for the log-generating function \(\loggen B z\) of
this process in terms of hypergeometric functions, see equation \eqref{eq:generating function}. This allows us to
compute limits and probabilities.
We recover the limit results by Kessler and Levin in \cite{kessler2013} and \cite{KessleronArXiV},
Durrett and Moseley \cite{durrettmoseley}, Mandelbrot \cite{mandelbrot1974} and and M\"ohle \cite{moehle2005}, and
extent them to a wider range of parameters and for included mutant-cell death, which is not treated in every case. 

To emphasize that the semi-deterministic model has many advantages over the fully stochastic one, we
compared simulation results of the fully stochastic model with the results of the LPSM-limit taken in Section
\ref{sec:limit behaviour}, see figures \ref{fig:Nmu1} and \ref{fig:Nmu10}. We find good agreement between
simulation results and the semi-stochastic theory already for relatively small values of \(N\).

We showed that the mutant distribution depends on the mutant extinction probability $q=\beta/\alpha$ in a nontrivial
way, and it cannot simply be scaled out of the formulas. In the large $\theta=N\mu$ limit, however, $q$ can be scaled
out of the mutant distribution for \(\gamma\in(0,2]\), hence only in this limit cell death has a trivial effect on the
mutant distribution.

From numerical work we found that the mutant distribution is unimodal, i.e.~it has only a single maximum. We gave a
simple condition for having no mutants being the most probable scenario. Otherwise, the most probable number of mutants
is positive, which we determined numerically and plotted for several parameter values.

We also determined the mean and the variance of the number of mutants, which are of course finite for finite wild-type
population size, $N$. When taking the large $N$ limit and using the appropriate scaling, we also obtain finite moments.
Conversely, in the large-population-small-mutation limit the moments become infinite for \(\gamma\in[0,1]\). This 
behavior, which caused some controversy in the literature (see \cite{zheng1999} for historical notes on the topic), is
explained by the power law tail of the mutant distribution in this limit. 

\begin{acknowledgement}

We thank the National Philantropic Trust \mbox{(FQEB Grant \#RFP-12-18)} for financial support.
\end{acknowledgement}

\appendix\normalsize

\section{Proof of recursion formulas}
\label{sec:rec}

To prove the recursion \eqref{eq:recursion_general} we first Taylor expand $\loggen{B}z$ around $z=0$ as
\begin{equation}
 \loggen{B}z
   = \sum_{n=0}^\infty q_n z^n
 \label{series_logG}
\end{equation}
Equating the coefficients of $z^n$ on both sides of 
\begin{equation}
 G(z)
   = \sum_{n=0}^\infty p_n z^n = e^{\sum_n^\infty q_n z^n}
\end{equation}
leads to the general recursion \eqref{eq:recursion_general} (stated as Lemma 3 in \cite{zheng1999}).

When Taylor expanding $\loggen{B}z $, we immediately obtain $q_0$ by replacing $\xi$ by $q=\xi|_{z=0}$ in
\eqref{eq:generating function}. To obtain higher order coefficients, we first calculate the coefficients $a_k$ of the
expansion of 
\begin{equation}
 \frac{1}{\gamma} \Fze{1,\gamma}{1+\gamma}{N^{-1/\gamma} \xi} = \sum_{k\ge 0} a_k z^k.
\end{equation}
With induction one can show that for $k\ge 1$
\begin{equation}
\begin{split}
& \frac{d^k}{dz^k}\Fze{1,\gamma}{1+\gamma}{N^{-1/\gamma} \xi} \\
& = k! \sum_{j=1}^k \binom{k-1}{j-1} \frac{\gamma (q-1)^j N^{-j/\gamma}}{(\gamma+j)(1-z)^{j+k}}
      \Fze{1+j,\gamma+j}{1+\gamma+j}{N^{-1/\gamma} \xi}.
\end{split}
\label{eq:derivatives first summand}
\end{equation}

Here the $k=1$ case can be easily checked, and the induction can be performed by differentiating the above expression
using \eqref{eq:derivative hypergeometric function}.  Now we take \eqref{eq:derivatives first
summand} at $z=0$ and simplify it by using \eqref{iden:to11} to find the coefficients
\begin{equation}
a_k = \sum_{j=1}^k \binom{k-1}{j-1} \frac{1}{j+\gamma} \left( \frac{1-q}{q-N^{1/\gamma}} \right)^j
         \Fze{1,\gamma}{1+\gamma+j}{N^{-1/\gamma}q}.
\end{equation}

For the second term in \eqref{eq:generating function} we need the $N=1$ case, where the above sum simplifies
\begin{equation}
\label{eq:ak_sum}
\begin{split}
a_k 
  & = \sum_{j=1}^k \binom{k-1}{j-1} \frac{( -1)^j}{j+\gamma}\Fze{1,\gamma}{1+\gamma+j}{q}\\
  & = \sum_{n\ge 0} q^n (\gamma)_n \sum_{j=1}^k \binom{k-1}{j-1} \frac{(-1)^j}{(j+\gamma)_{n+1}}.
\end{split}
\end{equation}
We can calculate the second sum by first generalizing it to a polynomial, re-indexing it, and realizing that we can
extend the summation to infinity 
\begin{equation}
\begin{split}
(\gamma)_{n}&\sum_{j=1}^k\binom{k-1}{j-1}\frac{(-1)^jz^{j-1}}{(j+\gamma)_{n+1}}\\
  & =-(\gamma)_n\sum_{j\ge 0}\binom{k-1}j \frac{(-1)^jz^j}{(j+\gamma+1)_{n+1}}\\
  & =-\frac{\gamma}{(\gamma+n)(\gamma+n+1)}\Fze{1-k,1+\gamma}{2+n+\gamma}z.
\label{eq:general identity}
\end{split}
\end{equation}
The parameters of the hypergeometric function were read out from the ratio of consecutive coefficients of the
series, as in \eqref{eq:Hyper_ratios}. 
Now taking the $z\nearrow1$ limit in \eqref{eq:general identity} and using the Chu-Vandermonde identity
\eqref{eq:Vander} yields
\begin{equation}
 (\gamma)_n \sum_{j=1}^k \binom{k-1}{j-1} \frac{(-1)^j}{(j+\gamma)_{n+1}} 
   = -\gamma \frac{(n+1)_{k-1}}{(n+\gamma)_{k+1}}.
\label{eq:binomial sum identity}
\end{equation}
We can substitute this expression into \eqref{eq:ak_sum} to obtain
\begin{equation}
\begin{split}
a_k 
  = -\gamma \sum_{n\ge 0}  q^n \frac{(n+1)_{k-1}}{(n+\gamma)_{k+1}}
  = -\frac{(k-1)!}{(\gamma+1)_{k}} \Fze{k,\gamma}{1+\gamma+k}{q},
\end{split}
\end{equation}
which immediately leads to \eqref{general_rec_coeffs}.

For the special case $\gamma=1$ we present an easier derivation.
We derive a recursion for the probabilities $p_n$ directly from the expression \eqref{b_one_gener} by Taylor
expanding around $z=0$ and using the binomial theorem to rewrite \(y^j\) as
\begin{equation}
 y^j 
   =(1-q)^{-j}\sum_{k=0}^j\binom jk(-q)^{j-k} z^k.
\end{equation}
Then
\begin{equation}
\begin{split}
\loggen{B}z
   & = \theta \frac{1-y}{y} \log (1-\phi y)\\
   & = -\theta \phi + \theta \sum_{j\ge1} \left( \frac{\phi^j}{j} - \frac{\phi^{j+1}}{j+1} \right) y^j\\
   & = -\theta\phi+\theta\sum_{j\ge 1} 
          \sum_{k=0}^j \binom{j}{k} \frac{(-q)^{j-k}}{(1-q)^j}\left(\frac{\phi^j}{j}
          -\frac{ \phi^{j+1}}{j+1} \right)z^k.
\end{split}
\label{neutral_logg_sum}
\end{equation}
By changing the order of summation, we read out the coefficients $q_k$ of $z^k$ as in \eqref{series_logG}
\begin{equation}
 q_k = 
 \begin{cases}
   -\theta\phi + \theta \sum_{j\ge 1} \left(\frac{-q}{1-q}\right)^j
        \left(
           \frac{\phi^j}{j}-\frac{\phi^{j+1}}{j+1}
        \right) & k=0,\\
   \frac{\theta \phi^k}{(1-q)^k} \sum_{j\ge k} \binom{j}{k}\left(\frac{-q}{1-q}\right)^{j-k}\left(\frac{\phi^{j-k}}{j}
         -\frac{\phi^{j-k+1}}{j+1} \right) & k\ge 1.
 \end{cases}
\end{equation}
For $q_0$ we have the series of a logarithm. For $k\ge 1$ we rewrite $q_k$ as two separate sums, and re-index the
summations from zero
\begin{equation}
 q_k 
   = \theta \left(\frac{\phi}{1-q}\right)^k \left[ \frac{1}{k} \sum_{j\ge 0} \binom{k-1+j}{k-1} x^j 
     + \phi \sum_{j\ge 0} \frac{(1+k)_j(1+k)_j}{(2+k)_j}\frac{x^j}{j!} \right]
\end{equation}
with $x=-q\phi/(1-q)$. The first sum is $(1-x)^{-k}$ and the second sum is a hypergeometric function, which can be
simplified using \eqref{iden:to11}, so finally we arrive at \eqref{eq:gamma1_rec}.

\section{Large $\gamma$-asymptotics of the resistance probability}\label{app:resistance}

We show the large \(\gamma\) approximations of the resistance probability \(p_0\) as given in
\eqref{eq:contourapproximation}. For this we utilize the following asymptotics for the hypergeometric function as
presented in \cite[15.12(iii)]{NIST:DLMF}. Adapted to our notation they state that for fixed \(a,b,c,z\in\mathds C\)
\begin{equation}
\Fze{a,b}{c+\gamma}z
  \sim \frac{\Gamma(c+\gamma)}{\Gamma(c-b+\gamma)}\sum_{k\geq0}r_k(z)(b)_k \gamma^{-k-b}
  \label{eq:asymptotic}
\end{equation}
for large \(\gamma\), where \(r_s(z)\) are the coefficients of the expansion of a specific function
\begin{equation}
 \sum_{k\geq 0}r_k(z)t^k=\left(\frac{e^t-1}t\right)^{b-1}e^{t(1-c)}(1-z+ze^{-t})^{-a}.
\end{equation}
We set for simplicity \(z=\frac{q}{q-1}\) and apply \eqref{eq:inversion hyperg 02} to the
hypergeometric function in \eqref{eq:contourapproximation}, then
\begin{equation*}
\begin{split}
   (1+\gamma)\Fze{1,\gamma}{2+\gamma}{q}^{-1}
   & = (1+\gamma)(1-q)\Fze{1,2}{2+\gamma}{z}^{-1}.
\end{split}
\label{eq:simplification hyperg}
\end{equation*}
We can now apply \eqref{eq:asymptotic} 
\begin{equation*}
 \Fze{1,2}{2+\gamma}{z}
  \sim \frac{\gamma+1}{\gamma}\sum_{k\geq0}r_k(z)\frac{(k+1)!}{\gamma^k}
\end{equation*}
where \(r_k(z)\) are the coefficient of the expansion of 
\begin{equation*}
 \frac{1-e^{-t}}{t}(1-z+ze^{-t})^{-1}=1+(z-1/2)t+(z^2-z+1/6)t^2+O(t^3).
\end{equation*}
Using this expansion in \eqref{eq:simplification hyperg} we arrive at \eqref{eq:contourapproximation}. Note that the
other similar results used throughout Section 5.1 on the resistence probability can be proven in the same
way and thus shall not be included here.

\section{Definitions and useful identities}\label{app:def and ids}

\emph{Euler's reflection formulas} for the Gamma function are
\begin{equation}
 \Gamma(1+z)\Gamma(1-z)=\frac{\pi z}{\sin(\pi z)},
 \label{eq:euler_reflection_B}
\end{equation}

\begin{equation}
 \Gamma(z+1)\Gamma(z-1)=\frac{z}{z-1} \Gamma(z)^2.
 \label{eq:euler_reflection_C}
\end{equation}

The Digamma function is defined as
$$
 \Psi(z) = \Psi_0(z) = \frac{d}{dz} \log \Gamma(z) = \frac{\Gamma'(z)}{\Gamma(z)}
$$
and has the expansion around $z=0$
\begin{equation}
\label{eq:series_Psi(0)}
 \Psi(z) = -\frac{1}{z} -E_C +\frac{\pi^2 z}{6} + O(z^2)
\end{equation}
and around $z=\infty$
\begin{equation}
\label{eq:series_Psi}
 \Psi(z) = \log z-\frac{1}{2 z}-\frac{1}{12 z^2}+O\left( \frac{1}{z^4}\right)
\end{equation}
Its generalization is the Polygamma function
$$
 \Psi_n(z) = \frac{d^{n+1}}{dz^{n+1}} \log \Gamma(z).
$$

The Pochhammer symbol is defined in terms of Gamma functions resp. as ascending factorial
\begin{equation}
 (a)_b =\frac{\Gamma(a+b)}{\Gamma(a)} = a(a+1)(a+2)\cdots(a+b-1).
\label{eq:definition pochhammer symbol}
\end{equation}

The series of the hypergeometric function is
\begin{equation}
 \Fze{a,b}{c}z = \sum_{n\geq 0}\frac{(a)_n (b)_n} {(c)_n}\frac{z^n}{n!}
\label{eq:definition gauss hypergeometric function}
\end{equation}
for any (complex) \(a,b,c\).
The coefficients of the power series $\sum_{n\geq 0} A_n z^n$ of the hypergeometric function \(\Fze{a,b}cz\) satisfy 
$A_0=1$ and 
\begin{equation}
\label{eq:Hyper_ratios}
 \frac{A_{n+1}}{A_n} = \frac{(n+a)(n+b)}{(n+c)(n+1)}.
\end{equation}
An alternative form of the hypergeometric function can be expressed as integral, if \(c-b>0\)
\begin{equation}
 \Fze{a,b}cz 
   = \frac{\Gamma(c)}{\Gamma(b)\Gamma(c-b)}\int_0^1\frac{t^{b-1}(1-t)^{c-b-1}}{(1-zt)^a}dt.
 \label{eq:hypergeo as integral}
\end{equation}

For $\Re(c-a-b)>0$ \emph{Gau\ss's hypergeometric Theorem} states
\begin{equation}
\label{eq:Vander}
  \Fze{a,b}{c}{1} = \frac{\Gamma(c)\Gamma(c-a-b)}{\Gamma(c-a)\Gamma(c-b)}
  = \frac{(c-a-b)_a}{(c-a)_a}
  = \frac{(c-a-b)_b}{(c-b)_b}
  = \frac{(c-b)_{-a}}{(c)_{-a}},
\end{equation}
which is also called the \emph{Chu-Vandermonde Identity} if $a$ is a negative integer.

A useful identity
\begin{equation}
 \Fze{a,b}cz=(1-z)^{-b}\Fze{c-a,b}c{\frac{z}{z-1}},
 \label{eq:inversion hyperg 02}
\end{equation}
when applied twice becomes
\begin{equation}
  \Fze{a,b}{c}{z} = (1-z)^{c-a-b} \Fze{c-a,c-b}{c}{z}.
  \label{iden:to11}
\end{equation}

\newcommand{\gammata}[4]{\frac{\Gamma(#1)\Gamma(#2)}{\Gamma(#3)\Gamma(#4)}}

Some inversion formulae (see \cite{NIST:DLMF})
\begin{equation}
\begin{split}
  \Fze{a,b}cz 
   & = \gammata c{b-a}b{c-a}
        (1-z)^{-a}\Fze{a,c-b}{a-b+1}{\frac1{1-z}}\\
   & + \gammata c{a-b}a{c-b}
        (1-z)^{-b}\Fze{b,c-a}{b-a+1}{\frac1{1-z}},
\end{split}
\label{eq:inversion hyperg 01}
\end{equation}

\begin{equation}
 \begin{split}
  \Fze{a,b}cz
   & = \gammata c{c-a-b}{c-b}{c-a}
        \Fze{a,b}{a+b-c+1}{1-z}\\
   & + \gammata c{a+b-c} ab
        z^{1-c}(1-z)^{c-a-b}\Fze{1-a,1-b}{c-a-b+1}{1-z}.
 \end{split}
\label{eq:hyperg trafo (1-z)}
\end{equation}

The derivative of the Hypergeometric function
\begin{equation}
\label{eq:derivative hypergeometric function}
 \frac{d}{dz}\Fze{a,b}{c}z=\frac{ab}{c}\Fze{a+1,b+1}{c+1}z.
\end{equation}

Occasionally we need the limit behavior for \(z\nearrow 1\). For general parameters the following formulas
hold:
If \(c=a+b\), then
\begin{equation}
 \lim_{z\nearrow 1}\frac{\Fze{a,b}{a+b}z}{-\log(1-z)}=\frac{\Gamma(a+b)}{\Gamma(a)\Gamma(b)}.
 \label{eq:z->1 c=a+b}
\end{equation}
If \(\Re(c-a-b)<0\), then
\begin{equation}
 \lim_{z\nearrow 1}\frac{\Fze{a,b}{a+b}z}{(1-z)^{c-a-b}}=\frac{\Gamma(c)\Gamma(a+b-c)}{\Gamma(a)\Gamma(b)}.
 \label{eq:z->1 re(c-a-b)<0}
\end{equation}
Note that the Chu-Vandermonde identity is also the limit case for \(\Re(c-a-b)>0\).

For a specific choice of the parameters the hypergeometric function can be expressed in simpler terms
\begin{equation}
 \Fze{1,1}2z=-\frac{\log(1-z)}{z}
 \label{eq:log and hypergeo}
\end{equation}

\begin{equation}
 \Fze{1,1}3z=\frac{2(z+(1-z)\log(1-z))}{z^2}
 \label{eq:log and hypergeo 2}
\end{equation}

For general \(z\) we can develop the hypergeometric function into
\begin{equation}
 \Fze{1,b}{c}z
   = 1+\frac{b}{c}z\Fze{1,b+1}{c+1}z
\label{eq:simplification}
\end{equation}
Note that this formula is not valid for \(a\neq1\).

\section{Theorems from Singularity Analysis}
\label{sec:singanal}

Theorem VI.1 \cite{flajolet_book}: (Standard function scale). Let $a$ be an arbitrary complex number in
$\mathbb{C}\setminus \mathbb{Z}_{\le0}$. The coefficient of $z^n$ for large $n$ has a full asymptotic expansion in
descending powers of $n$ 
\begin{equation}
\label{eq:fla1}
\begin{split}
 [z^n] (1-z)^{-a}  
 &\sim \frac{n^{a-1}}{\Gamma(a)} \left( 1+\sum_{k\ge1} \frac{e_k}{n^k} \right)\\
 &\sim \frac{n^{a-1}}{\Gamma(a)} 
 \left( 1 + \frac{a(a-1)}{2n} + \frac{a(a-1)(a-2)(3a-1)}{24 n^2} + \dots \right)
\end{split}
\end{equation}
where $e_k$ is a polynomial in $a$ of degree $2k$ and specified in \cite{flajolet_book}.

Theorem VI.2 \cite{flajolet_book}: (Standard function scale, logarithms). Let $a$ be an arbitrary complex number in
$\mathbb{C}\setminus \mathbb{Z}_{\le0}$ and $b\in\mathbb{C}$. The coefficient of $z^n$ for large $n$ has a full 
asymptotic expansion in descending powers of $n$
\begin{equation}
\label{eq:fla2}
 [z^n] (1-z)^{-a} \left( \frac{L}{z} \right)^b \sim \frac{n^{a-1}}{\Gamma(a)} (\log n)^b 
 \left[ 1+ \frac{C_1}{\log n} + \frac{C_2}{\log^2 n} + \dots \right],
\end{equation}
where $C_k = \binom{b}{k}\Gamma(a) \frac{d^k}{ds^k} \frac{1}{\Gamma(s)}\big|_{s=a}$. (Note that there is an erroneous
$(-1)^k$ factor in the online version of the book.)

For $b=1$ we only have two terms
\begin{equation}
\label{eq:fla2b}
 [z^n] (1-z)^{-a} \frac{L}{z}
 \sim \frac{n^{a-1}}{\Gamma(a)}  
  \left[ \log n - \Psi(a) \right],
\end{equation}
where $\Psi(a)=\Gamma'(a)/\Gamma(a)$ is the Digamma function.

For $a\in \mathbb{Z}_{\le 0}$ and $b\in \mathbb{Z}_{\ge 0}$
\begin{equation}
\label{eq:fla_ints}
 [z^n] (1-z)^{-a} L^b \sim 
 n^{a-1} \left[ F_0(\log n) + \frac{F_1(\log n)}{n} + \dots \right],
\end{equation}
 where the degree of $F_n$ is $k-1$. The polynomials are given by the following Frobenius-Jungen method. For $b\in
\mathbb{Z}_{\ge 0}$ and arbitrary $a\in\mathbb{C}$ 
\begin{equation}
\label{eq:fla_bint}
 [z^n] (1-z)^{-a} L^b = 
 \frac{d^b}{da^b} \binom{n+a-1}{n}.
\end{equation}
Note also that \eqref{eq:fla_bint} is exact apart from finitely many terms. 
For example for $b=1$ and 2, it takes the form
\begin{equation}
\label{eq:fla_b12agen}
\begin{split}
 [z^n] (1-z)^{-a} L &= \binom{n+a-1}{n} h_n(a),\\
 [z^n] (1-z)^{-a} L^2 &= \binom{n+a-1}{n} [h'_n(a)-h_n^2(a)],\\
\end{split} 
\end{equation}
where
\begin{equation}
\label{eq:fla_h}
\begin{split}
h_n(a)&=\Psi(a+n)+\Psi(n),\\ 
h'_n(a)&=\Psi_1(a+n)+\Psi_1(n).\\ 
\end{split} 
\end{equation}
For negative integer $a$ consider the expression as a limit.
For example for $b=1$ it simplifies to
\begin{equation}
\label{eq:fla_b1aint}
  [z^n] (1-z)^k L  = (-1)^k \frac{k!}{n(n-1)\cdots (n-k)} 
\end{equation}
and for $b=2$ 
\begin{equation}
\label{eq:fla_b2aint}
\begin{split}
  [z^n] (1-z)^k L^2  &= 2 (-1)^k k! \frac{C_E+\Psi(n-k)-H_k}{n(n-1)\cdots(n-k)},
\end{split}  
\end{equation}
where $H_n=\sum_{k=1}^n1/n$ is the $n$-th harmonic number.

\bibliographystyle{ieeetr}
\bibliography{bibliography}

\end{document}